\documentclass[11pt]{amsart}
\usepackage{amssymb,latexsym}
\usepackage{pstcol,pstricks,color}

 \numberwithin{equation}{section}

\newtheorem{theorem}{Theorem}[section]
\newtheorem{proposition}[theorem]{Proposition}

\newtheorem{remark}[theorem]{Remark}

\newtheorem{lemma}[theorem]{Lemma}

\title{On the number of combinations without certain separations}

\begin{document}
\maketitle
\begin{center}
Toufik Mansour$^\dag$ and Yidong Sun$^\ddag$

$^\dag$Department of Mathematics, University of Haifa, 31905 Haifa,
Israel\\
$^\ddag$Department of Mathematics, Dalian Maritime University, 116026 Dalian, P.R. China\\[5pt]

{\it $^\dag$toufik@math.haifa.ac.il, $^\ddag$sydmath@yahoo.com.cn}
\end{center}\vskip0.5cm

\subsection*{Abstract}
In this paper we enumerate the number of ways of selecting $k$
objects from $n$ objects arrayed in a line such that no two selected
ones are separated by $m-1,2m-1,\cdots,pm-1$ objects and provide
three different formulas when $m,p\geq 1$ and $n\geq pm(k-1)$. Also,
we prove that the number of ways of selecting $k$ objects from $n$
objects arrayed in a circle such that no two selected ones are
separated by $m-1,2m-1,\cdots,pm-1$ objects is given by
$\frac{n}{n-pk}\binom{n-pk}{k}$, where $m,p\geq 1$ and $n\geq
mpk+1$.
\medskip

{\bf Keywords}: Composition, Cauchy's residue theorem, $N$-separate

\noindent {\sc 2000 Mathematics Subject Classification}: Primary
05A05, 05A15

\section{Introduction}
In 1943, Kaplansky \cite{Kaplansky} published a recursive derivation
of the number of combinations of $n$ objects taken $k$ at a time
without two selected ones being consecutive (see also Comtet
\cite{Comtet}, Riordan \cite{Riordan} and Ryser \cite{Ryser}). In
1981, Konvalina \cite{Konvalina} derived the number of combinations
of $n$ objects taken $k$ at a time without two selected ones having
unit separation, i.e., having exactly one object between them.

Let $[n]\ ({\rm resp.}\ [\overline{n}])$ be the set of $n$ objects
$x_1,x_2,\cdots,x_n$ arrayed in a line (resp.\ circle). Given a
subset $N$ of the set $\mathbb{N}$ of nonnegative integers, a subset
$A$ of $[n]$ or $[\overline{n}]$ will be called {\em $N$-separate}
if any two objects in $A$ have exactly $j$ objects between them,
then $j\in N$. Let $N_{m}^p=\mathbb{N}-\{m-1,2m-1,\cdots,pm-1\}$ for
any integers $m,p\geq 1$, define $\mathcal{H}_{p,n}^{m,k}\ ({\rm
resp.}\ \mathcal{G}_{p,n}^{m,k}$) to be the number of
$N_{m}^p$-separate $k$-subsets of $[n]\ ({\rm resp.}\
[\overline{n}]$). Thus, by our notation, Konvalina \cite{Konvalina}
considered the special case $N_{2}^1$-separation of $[n]$ and
$[\overline{n}]$, Kaplansky \cite{Kaplansky} discussed the special
case $N_{1}^p$-separation of $[n]$ and $[\overline{n}]$, and
obtained that
\begin{eqnarray}\label{eqn 1.1}
\mathcal{H}_{p,n}^{1,k}=\binom{n-p(k-1)}{k}\quad\mbox{ and }\quad
\mathcal{G}_{p,n}^{1,k}=\frac{n}{n-pk}\binom{n-pk}{k}.
\end{eqnarray}
In this paper, by combinatorial analysis together with algebraic
method, we extend the above results to the general case of $m$.

\section{Some preliminary remarks}
Let $n=rm+\ell$ with $1\leq \ell\leq m$, and let $A_1,\ldots, A_m$
be a partition of $[n]=\{x_1,x_2,\cdots,x_n\}$ into $m$ subsets
defined
$$\begin{array}{lll}
A_i&=\{x_i,x_{m+i},\cdots,x_{rm+i}\},& 1\leq i\leq \ell, \\
A_i&=\{x_i,x_{m+i},\cdots,x_{(r-1)m+i}\},& \ell+1\leq i\leq m,
\end{array}$$
then put them in an array,
\begin{equation*}
\begin{array}{ccccc}
x_1 & x_{m+1} & \cdots & x_{(r-1)m+1} & x_{rm+1} \\
\vdots & \vdots &  \vdots  & \vdots & \vdots \\
x_l & x_{m+l} & \cdots & x_{(r-1)m+l} & x_{rm+l}\\
x_{l+1} & x_{m+l+1} & \cdots & x_{(r-1)m+l+1} & \empty  \\
\vdots & \vdots &  \vdots  & \vdots &  \empty   \\
x_m & x_{2m} & \cdots & x_{(r-1)m+m} & \empty
\end{array}
\end{equation*}

For any $k$-subset $B$ of $[n]$, define $B_i=B\cap A_i$. Note that,
in the line case, $B$ is $N_{m}^p$-separate if and only if each
$B_i$ is $N_{1}^p$-separate. From this critical observation together
with (\ref{eqn 1.1}), we can obtain the following result.
\begin{proposition}\label{prop 2.1}
For any integers $p,m\geq 1$ and $n,k\geq 0$,
\begin{eqnarray}\label{eqn
2.1}
\mathcal{H}_{p,n}^{m,k}&=&\sum_{\sigma_1(k,m)}\,\prod_{i=1}^m\binom{|A_i|-p(k_i-1)}{k_i},
\end{eqnarray}
where $|A_i|$ is the cardinality of the set $A_i$, and
$\sigma_1(k,m)$ denotes the all nonnegative integer solutions of
$k_1+k_2+\cdots+k_m=k$ such that $k_i\leq 1+\frac{|A_i|}{p}$ for
$i=1,2,\ldots,m$.
\end{proposition}

In the next section, we can find the explicit formula for
$\mathcal{H}_{p,n}^{m,k}$, and show that when $n$ is large enough
$(n\geq mp(k-1)\ {\rm here})$, then $\mathcal{H}_{p,n}^{m,k}$ is
independent of the composition of $n$, i.e.,
$|A_1|+|A_2|+\cdots+|A_m|=n$. However, in the circle case, the above
decomposition does not work, for example, when $n=5,p=1,m=k=2$, then
$[\overline{5}]=\{x_1,x_2,\cdots,x_5\}$ has five $N_{2}^1$-separate
$2$-subsets, which are $\{x_1,x_2\},
\{x_2,x_3\},\{x_3,x_4\},\{x_4,x_5\},\{x_5,x_1\}$, while
$\{x_5,x_1\}\cap \{x_1,x_3,x_5\}=\{x_5,x_1\}$ is not an
$N_{1}^1$-separate $2$-subset of $\{x_1,x_3,x_5\}$. In spite of
this, we can derive a recurrence relation between
$\mathcal{H}_{p,n}^{m,k}$ and $\mathcal{G}_{p,n}^{m,k}$ for $n\geq
mpk+1$.

Given a $N_{m}^p$-separate $k$-subset $B$ of $[\overline{n}]$, for
some $0\leq j\leq m$, there exist $j$ elements of $B$, say
$x_{i_1},x_{i_2},\cdots,x_{i_j}$, lying in the subset
$\{x_1,x_2,\cdots,x_{mp}\}$, in other words, each of which is
respectively one of the first $p$ objects of
$A_{\ell_1},A_{\ell_2},\cdots,A_{\ell_j}$, then there are
$\binom{m}{j}p^j$ ways to do this. Now delete the related $j(2p+1)$
objects of $[\overline{n}]$, and delete the remainder $p(m-j)$
elements of $\{x_1,x_2,\cdots,x_{mp}\}$, then we get $m$ object sets
$A'_1,A'_2,\cdots,A'_m$ in which all elements are arrayed in a line
and there are totally $n-p(m-j)-j(2p+1)=n-pm-pj-j$ elements. Note
that the condition $n\geq mpk+1$ leads to $n-pm-pj-j\geq mp(k-j-1)$,
which makes the restricted inequality condition of \eqref{eqn 2.1}
in Proposition \ref{prop 2.1} redundant. Then there are
$\mathcal{H}_{p,n-pm-pj-j}^{m,k-j}$ ways to select the other $k-j$
objects from $A'_1,A'_2,\cdots,A'_m$. Hence, we have
\begin{proposition}\label{prop 2.2}
For any integers $p,m\geq 1$, $n,k\geq 0$ and $n\geq mpk+1$,
\begin{eqnarray}\label{eqn 2.2}
\mathcal{G}_{p,n}^{m,k}&=&\sum_{j\geq
0}\binom{m}{j}p^j\mathcal{H}_{p,n-pm-pj-j}^{m,k-j}.
\end{eqnarray}
\end{proposition}

Clearly, we can easily compute special values for
$\mathcal{H}_{p,n}^{m,k}$ and $\mathcal{G}_{p,n}^{m,k}$, that is,
\begin{itemize}
\item $\mathcal{H}_{p,n}^{m,k}=\mathcal{G}_{p,n}^{m,k}=0$ for $n< k$;

\item $\mathcal{H}_{p,n}^{m,0}=\mathcal{G}_{p,n}^{m,0}=1$;

\item $\mathcal{H}_{p,n}^{m,1}=\mathcal{G}_{p,n}^{m,1}=n$ for $n\geq
1$;

\item $\mathcal{H}_{p,n+k}^{m,k}=0$ for $im+1\leq k\leq (i+1)m$,
$0\leq n< imp$ and $i\geq 1$;

\item $\mathcal{G}_{p,n+k}^{m,k}=0$ for $im+1\leq k\leq (i+1)m$,
$0\leq n< (i+1)mp$ and $i\geq 1$.
\end{itemize}
Define $\mathcal{H}_{p,n}^{m,k}=\mathcal{G}_{p,n}^{m,k}=0$ for $k<0$
or $n<0$.

\section{Main result}
In order to give explicit formulas for $\mathcal{H}_{p,n}^{m,k}$ and
$\mathcal{G}_{p,n}^{m,k}$, we need the following critical lemma.

\begin{lemma}\label{lemma 3.1}
Let $\lambda_1,\lambda_2,\cdots,\lambda_m,\mu$ be any $m+1$ complex
numbers and $\lambda=\lambda_1+\lambda_2+\cdots+\lambda_m$. Define
$$\begin{array}{ll}
\Omega_{\mu,\lambda}^{m,k}(\lambda_1,\lambda_2,\cdots,\lambda_m)&=\sum\limits_{\sigma(k,m)}\prod\limits_{i=1}^m\binom{\lambda_i+\mu
k_i}{k_i},\\
\Phi_{\mu,\lambda}^{m,k}(\lambda_1,\lambda_2,\cdots,\lambda_m)&=
\sum\limits_{\sigma(k,m)}\prod\limits_{i=1}^m\frac{\lambda_i}{\lambda_i+\mu
k_i}\binom{\lambda_i+\mu k_i}{k_i}, \end{array}$$ where
$\sigma(k,m)$ denotes the all nonnegative integer solutions of
$k_1+k_2+\cdots+k_m=k$.

Then for all $m\geq 1$ and $n,k\geq 0$,
{\small\begin{eqnarray}
\Omega_{\mu,\lambda}^{m,k}(\lambda_1,\lambda_2,\cdots,\lambda_m)
&=&\sum_{j\geq 0}\binom{m+j-2}{j}\binom{\lambda+\mu k+m-1}{k-j}(\mu-1)^j, \label{eqn 3.1}\\
&=&\sum_{j\geq 0}\binom{\lambda+(\mu-1)k+j}{j}\binom{\lambda+\mu k+m-1}{k-j}(1-\mu)^j\mu^{k-j},\label{eqn 3.2}\\
&=&\sum_{j\geq 0}\frac{\lambda+\mu(m+j)}{k}\binom{m+j-1}{j}\binom{\lambda+\mu k+m-1}{k-j}(\mu-1)^j,\label{eqn 3.3}\\
\Phi_{\mu,\lambda}^{m,k}(\lambda_1,\lambda_2,\cdots,\lambda_m) &=&
\frac{\lambda}{\lambda+\mu k}\binom{\lambda+\mu k}{k}.\label{eqn
3.4}\end{eqnarray}}
\end{lemma}
\begin{proof}
First we recall the definition of the residue of a function. Let
$z_0$ be any isolated singular point of a function $f$. Then there
is a Laurent series $f(z)=\sum_{j=-\infty}^\infty a_j(z-z_0)^j$
valid for $0<|z-z_0|<R$, for some positive $R$. The coefficient
$a_{-1}$ of $(z-z_0)^{-1}$ is called the residue of $f$ at $z_0$,
and is usually written $\underset{z=z_0}{Res}f$ (for computing and
properties of the residue see for example \cite{Egorychev,Henrici}).
For simplicity, we write $\underset{z}{Res}f$ instead
$\underset{z=0}{Res}f$.

Note that the generalized binomial coefficient $\binom{\lambda}{k}$
has an integral representation,
\begin{eqnarray*}
\binom{\lambda}{k}=\underset{x}{Res}\frac{(1+x)^\lambda}{x^{k+1}},
\end{eqnarray*}
which yields that
\begin{eqnarray}\label{eqn residue}
\frac{\lambda}{\lambda+\mu k}\binom{\lambda+\mu
k}{k}=\underset{x}{Res}\frac{(1+x)^{\lambda+\mu
k-1}(1-(\mu-1)x)}{x^{k+1}}.
\end{eqnarray}
Then we have
\begin{eqnarray*}
\Omega_{\mu,\lambda}^{m,k}(\lambda_1,\lambda_2,\cdots,\lambda_m)&=&\underset{x}{Res}\frac{1}{x^{k+1}}
\prod_{i=1}^m\sum_{k_i\geq
0}x^{k_i}\underset{y_i}{Res}\frac{(1+y_i)^{\lambda_i+\mu
k_i}}{y_i^{k_i+1}},\\
&=&\underset{x}{Res}\left\{\prod_{i=1}^{m}\frac{(1+y_i)^{\lambda_i+1}}
{1-(\mu-1)y_i}\Big|_{y_i=x(1+y_i)^{\mu}}\right\}x^{-k-1},\\
&=&\underset{x}{Res}\frac{(1+\varphi(x))^{\lambda+m}}{(1-(\mu-1)\varphi(x))^m}x^{-k-1},
\end{eqnarray*}
where $\varphi(x)=x(1+\varphi(x))^{\mu}$. Using the Lagrange
inversion formula for $k\geq1$ with replacing $x$ by
$\frac{y}{(1+y)^{\mu}}$, we get that
\begin{eqnarray*}
\Omega_{\mu,\lambda}^{m,k}(\lambda_1,\lambda_2,\cdots,\lambda_m)&=&\sum_{j\geq 0}\frac{\lambda+\mu(m+j)}{k}\binom{m+j-1}{j}\binom{\lambda+\mu k+m-1}{k-j}(\mu-1)^j, \\
&=&\underset{y}{Res}\frac{(1+y)^{\lambda+\mu
k+m-1}}{(1-(\mu-1)y)^{m-1}}y^{-k-1}, \\
&=&\sum_{j\geq 0}\binom{m+j-2}{j}\binom{\lambda+\mu k+m-1}{k-j}(\mu-1)^j,\\
&=&\underset{y}{Res}(1-(\mu-1)y)^{\lambda+\mu k}(1+\frac{\mu
y}{1-(\mu-1)y})^{\lambda+\mu k+m-1}y^{-k-1},\\
&=&\sum_{j\geq 0}\binom{\lambda+(\mu-1)k+j}{j}\binom{\lambda+\mu k+m-1}{k-j}(1-\mu)^j\mu^{k-j}.\\
\end{eqnarray*}
Similarly, we also have
\begin{eqnarray*}
\Phi_{\mu,\lambda}^{m,k}(\lambda_1,\lambda_2,\cdots,\lambda_m)&=&\underset{x}{Res}\frac{1}{x^{k+1}}
\prod_{i=1}^m\sum_{k_i\geq
0}x^{k_i}\underset{y_i}{Res}\frac{(1+y_i)^{\lambda_i+\mu
k_i-1}(1-(\mu-1)y_i)}{y_i^{k_i+1}},\\
&=&\underset{x}{Res}\frac{(1+x)^{\lambda+\mu
k-1}(1-(\mu-1)x)}{x^{k+1}},\\
&=&\frac{\lambda}{\lambda+\mu k}\binom{\lambda+\mu k}{k}.
\end{eqnarray*}
This completes the proof.
\end{proof}

\begin{remark}
Note that Hwang and Wei \cite{Hwang} considered the special case
$$\Omega_{-1,n+m}^{m,k}(n_1+1,\cdots,n_m+1)=\sum_{\sigma(k,m)}\prod_{i=1}^m\binom{n_i+1-k_i}{k_i},$$
with $n=n_1+n_2+\cdots+n_m$ and obtained its another expression by
recurrence relation,
\begin{eqnarray*}
\Omega_{-1,n+m}^{m,k}(n_1+1,\cdots,n_m+1)&=&\sum_{j\geq
0}\binom{m+j-2}{j}\binom{n+1-k-2j}{k-2j}, \nonumber
\end{eqnarray*}
which can be derived easily from the proof of Lemma \ref{lemma 3.1}
if one notices that
$$\begin{array}{l}
\Omega_{-1,n+m}^{m,k}(n_1+1,\cdots,n_m+1)
=\underset{y}{Res}\frac{(1+y)^{n+2m-k-1}}{(1+2y)^{m-1}}y^{-k-1}
=\underset{y}{Res}\frac{(1+y)^{n-k+1}}{(1-\frac{y^2}{(1+y)^2})^{m-1}}y^{-k-1}.
\end{array}$$
Also, the equation \eqref{eqn 3.4} is a generalization of Gould's
identity \cite{Blackwell,Comtet}, that is,
\begin{eqnarray*}
\sum_{k=0}^n\frac{a}{a+ck}\binom{a+ck}{k}\frac{b}{b+c(n-k)}\binom{b+c(n-k)}{n-k}
&=&\frac{a+b}{a+b+ck}\binom{a+b+ck}{k}.
\end{eqnarray*}
Then (\ref{eqn 3.4}) can be proved again by repeatedly using Gould's
identity.
\end{remark}

Notice that when $n\geq pm(k-1)$ in (\ref{eqn 2.1}), then the
inequality condition for $\sigma_1(k,m)$ (i.e., $k_i\leq 1+
\frac{|A_i|}{p}$) is redundant. Hence, setting $\lambda_i=|A_i|+p$,
$\mu=-p$ in (\ref{eqn 3.1})--(\ref{eqn 3.3}), and combining with
Proposition \ref{prop 2.1}, we obtain our main result.

\begin{theorem} Let $p,m,k\geq 1$ be any integers. For $n\geq
pm(k-1)$,
\begin{eqnarray*}
\mathcal{H}_{p,n}^{m,k}&=&\sum_{j\geq 0}\binom{m+j-2}{j}\binom{n+mp+m-pk-1}{k-j}(-p-1)^j,\\
&=&\sum_{j\geq 0}\binom{n+mp-(p+1)k+j}{j}\binom{n+mp+m-pk-1}{k-j}(p+1)^j(-p)^{k-j},\\
&=&\sum_{j\geq
0}\frac{n-pj}{k}\binom{m+j-1}{j}\binom{n+mp+m-pk-1}{k-j}(-p-1)^j,
\end{eqnarray*}
and for $n\geq mpk+1$,
\begin{eqnarray}\label{eqn 3.5}
\mathcal{G}_{p,n}^{m,k}&=&\frac{n}{n-pk}\binom{n-pk}{k}.
\end{eqnarray}
\end{theorem}
\begin{proof}
It just needs to prove (\ref{eqn 3.5}). For $n\geq mpk+1$, by
(\ref{eqn 2.2}), we have
\begin{eqnarray*}
\mathcal{G}_{p,n}^{m,k}&=&\sum_{j\geq
0}\binom{m}{j}p^j\mathcal{H}_{p,n-pm-pj-j}^{m,k-j}\\
&=&\sum_{j\geq
0}\binom{m}{j}p^j\underset{y}{Res}\frac{(1+y)^{n-p(k-j)+m-1-pj-j}}{(1+(p+1)y)^{m-1}}y^{-(k-j)-1} \\
&=&\underset{y}{Res}\frac{(1+y)^{n-pk+m-1}}{(1+(p+1)y)^{m-1}}y^{-k-1}\sum_{j\geq
0}\binom{m}{j}p^j\left\{\frac{y}{1+y}\right\}^j \\
&=&\underset{y}{Res}\frac{(1+y)^{n-pk-1}(1+(p+1)y)}{y^{k+1}}\\
&=&\frac{n}{n-pk}\binom{n-pk}{k},
\end{eqnarray*}
which follows by (\ref{eqn residue}).
\end{proof}

The formulas (\ref{eqn 1.1}) and (\ref{eqn 3.5}) motivate the
following
\begin{theorem}
For any integers $p,m,n,k\geq 1$, if $n\geq mpk+1$, then there
exists a bijection between the set of $N_{m}^p$-separate $k$-subsets
of $[\overline{n}]$ and the set of $N_{1}^p$-separate $k$-subsets of
$[\overline{n}]$.
\end{theorem}

We fail to produce such a bijection, and find it remains a
challenging open question.

Now, we give several recurrence relations that the sequences
$\mathcal{H}_{p,n}^{m,k}$ and $\mathcal{G}_{p,n}^{m,k}$ satisfy.

\begin{theorem} Let $p,m,k\geq 1$ be any integers. For $n\geq pm(k-1)$,
\begin{eqnarray}\label{eqn 4.1}
\mathcal{H}_{p,n}^{m,k}&=&\mathcal{H}_{p,n-1}^{m,k}+\mathcal{H}_{p,n-p-1}^{m,k-1},
\end{eqnarray}
and for $n\geq m(pk+1)$,
\begin{eqnarray}
\mathcal{G}_{p,n}^{m,k}&=&\mathcal{G}_{p,n-1}^{m,k}+\mathcal{G}_{p,n-p}^{m,k-1},\label{eqn 4.2}\\
\mathcal{G}_{p,n}^{m,k}&=&\sum_{j\geq
0}(-1)^j\binom{m}{j}p^j(p+1)^{m-j}\mathcal{H}_{p,n-pm-j}^{m,k},\label{eqn
4.4}
\end{eqnarray}
and for $n\geq mp(k-1)$,
\begin{eqnarray}
\mathcal{H}_{p,n}^{m,k}&=&\sum_{j\geq
0}(-1)^j\binom{m+j-1}{j}p^j\mathcal{G}_{p,n+pm-pj-j}^{m,k-j}.\label{eqn
4.5}
\end{eqnarray}
\end{theorem}
\begin{proof}
To prove \eqref{eqn 4.1}, let us consider $N_{m}^p$-separate
$k$-subsets from $[n]$ which either contain the first object $x_1$
or do not. In the later case, the number of such subsets is
enumerated by $\mathcal{H}_{p,n-1}^{m,k}$. In the former case, the
subsets does not contain the objects
$x_{m+1},x_{2m+1},\cdots,x_{pm+1}$ of the set $A_1$ as sefined in
Section 2, note that the condition $n\geq mp(k-1)$ makes the
restricted inequality condition of \eqref{eqn 2.1} in Proposition
\ref{prop 2.1} redundant, so such subsets are counted by
$\mathcal{H}_{p,n-p-1}^{m,k-1}$. Hence, \eqref{eqn 4.1} holds.

Using simple algebraic calculations we obtain that \eqref{eqn 4.2}
holds.

Note that if $n\geq m(pk+1)$, there holds
\begin{eqnarray*}
\mathcal{G}_{p,n}^{m,k}&=&\underset{x}{Res}\frac{(1+x)^{n-pk-1}(1+(p+1)x)}{x^{k+1}} \\
&=&\underset{x}{Res}\frac{(1+x)^{n+m-pk-1}}{(1+(p+1)x)^{m-1}}\Big(p+1-\frac{p}{1+x}\Big)^mx^{-k-1}\\
&=&\sum_{j\geq
0}(-1)^j\binom{m}{j}p^j(p+1)^{m-j}\underset{x}{Res}\frac{(1+x)^{n+m-pk-j-1}}{(1+(p+1)x)^{m-1}}x^{-k-1}\\
&=&\sum_{j\geq
0}(-1)^j\binom{m}{j}p^j(p+1)^{m-j}\mathcal{H}_{p,n-pm-j}^{m,k},
\end{eqnarray*}
and if $n\geq mp(k-1)$, there holds
\begin{eqnarray*}
\mathcal{H}_{p,n}^{m,k}&=&\underset{x}{Res}\frac{(1+x)^{n+pm+m-pk-1}}{(1+(p+1)x)^{m-1}}{x^{-k-1}} \\
&=&\underset{x}{Res}\frac{(1+x)^{n+pm-pk-1}(1+(p+1)x)}{x^{k+1}}\Big({1+\frac{px}{1+x}}\Big)^{-m}\\
&=&\sum_{j\geq
0}(-1)^j\binom{m+j-1}{j}p^j\underset{x}{Res}\frac{(1+x)^{n+pm-pj-j-p(k-j)-1}(1+(p+1)x)}{x^{k-j+1}}\\
&=&\sum_{j\geq
0}(-1)^j\binom{m+j-1}{j}p^j\mathcal{G}_{p,n+pm-pj-j}^{m,k-j},
\end{eqnarray*}
which prove \eqref{eqn 4.4} and \eqref{eqn 4.5}.
\end{proof}

The above theorem suggests that there should exist combinatorial
proofs for \eqref{eqn 4.2}, \eqref{eqn 4.4} and \eqref{eqn 4.5}.
However, we fail to produce such proofs, and find them remain
challenging open questions.

\vskip0.5cm

\subsection*{Acknowledgment} We are grateful to Victor
J. W. Guo for helpful discussions. Thanks also to the referee for
valuable suggestions.

\vskip0.5cm

%==============================================================================================================

\end{document}